# Semi-primitive roots and irreducible quadratic forms

24th July 2024

**WOLF Marc,** https://orcid.org/0000-0002-6518-9882

Email: marc.wolf3@wanadoo.fr

**WOLF François,** https://orcid.org/0000-0002-3330-6087

Email: francois.wolf@dbmail.com

**Abstract**

Modulo a prime number, we define semi-primitive roots as the square of primitive roots. We present a method for calculating primitive roots from quadratic residues, including semi-primitive roots. We then present progressions that generate primitive and semi-primitive roots, and deduce an algorithm to obtain the full set of primitive roots without any $gcd$ calculation.

Next, we present a method for determining irreducible quadratic forms with arbitrarily large conjectured asymptotic density of primes (after Shanks, [1][2]). To this end, we propose an algorithm for calculating the square root modulo p, based on the Tonelli-Shanks algorithm [4].

CONTENTS

# 1. Introduction

## 1.1. Primitive root modulo N

Quadratic residues and non-residues, as well as primitive roots, have given rise to an abundant literature in mathematics (group theory with Lagrange's and Fermat's theorems, ring and field theory, etc.). They have many applications, including cryptography, primality tests and integer factorisation.

We recall that the group of units of the ring $\mathbb{Z}/N\mathbb{Z}$, of order $\varphi(N)$ (Euler's totient of $N$), is a cyclic group if and only if $N = 2$, $N = 4$, $N = p^k$ or $N = 2p^k$ with $p \in \mathbb{P}\setminus\{2\}$ and $k \in \mathbb{N}^*$. For such $N$, a <u>primitive root modulo $N$</u> is a generator of this group. There exist $\varphi(\varphi(N))$ such integers modulo $N$, whose powers generate every element of $(\mathbb{Z}/N\mathbb{Z})^\times$. The usual method to compute one or more primitive roots requires the prime decomposition of $\varphi(N)$. In this article, we will focus on the case $N = p \in \mathbb{P} \setminus \{2\}$.

<u>In the first section</u>, we define semi-primitive roots and study their relationship with primitive roots modulo $p$. Using these results, we present a test on quadratic residues to determine primitive roots. Then, using recursive sequences of primitive and semi-primitive roots, we describe an algorithm determining all primitive roots without any $gcd$ computation.

In [2], we studied the density of primes of the form $X^2 + c$, for a fixed $c \in \mathbb{N}^*$, based on Shanks' conjecture [1], which we empirically corroborated. <u>In the second section of this article</u>, we will continue this study and propose new results generalised to irreducible quadratic forms.

## 1.2. Basic properties, definitions and notations

**Definition 1.1**:

a) We let $\mathcal{D}(x) := \{d \in \mathbb{N}^* | d|x\}$ the set of divisors of an integer $x$, for $E \subset \mathbb{Z}$ we let $\mathcal{D}(E) := \bigcup_{x \in E} \mathcal{D}(x)$ the set of divisors of at least one element of $E$. We also define $\mathcal{D}_p(E) = \mathcal{D}(E) \cap \mathbb{P}$ the set of prime divisors of $E$.
b) For $(x, y) \in \mathbb{Z}^* \times \mathbb{N}^*$ we let $\mathcal{P}_{x,y} := \{p \in \mathbb{P} \mid p \equiv x\ [y]\}$.
c) Let $p \in \mathbb{P} \setminus \{2\}$. A <u>semi-primitive root modulo $p$</u> is defined as the square of a primitive root in $(\mathbb{Z}/p\mathbb{Z})^\times = \mathbb{F}_p^*$.

**Property 1.1**: $g \in \mathbb{F}_p^*$ is a semi-primitive root if and only if the order of $\langle g \rangle = \{g^k, k \in \mathbb{Z}\}$ is $\frac{p-1}{2}$.

*Proof*: If $g$ is a semi-primitive root, there exists a primitive root $h$ such as $g = h^2$. Since $p - 1$ is even and the order of $\langle h \rangle$ is $p - 1$, we deduce that the order of $\langle g \rangle$ is $\frac{p-1}{2}$. Conversely, if the order of $\langle g \rangle$ is $\frac{p-1}{2}$ and if $h$ is a primitive root, let $k$ such that $g = h^k$. We know that the order of $\langle g \rangle$ equal to $\frac{p-1}{\gcd(k,p-1)}$ thus $\gcd(k, p - 1) = 2$, i.e. we can write $k = 2u$ with $u$ and $p - 1$ coprime. We conclude that $h^u$ is a primitive root and $g = (h^u)^2$.

**Definition 1.2**: Let $p \in \mathbb{P} \setminus \{2\}$. $\mathcal{G}_{Z,p}$ (respectively $\mathcal{G}_{S,p}$, $\mathcal{G}_{QR,p}$, $\mathcal{G}_{QNR,p}$) is the set of primitive roots (respectively semi-primitive roots, quadratic residues, quadratic non-residues) modulo $p$. When the value of p is unambiguous, this set is simply noted $\mathcal{G}_Z$ (respectively $\mathcal{G}_S$, $\mathcal{G}_{QR}$, $\mathcal{G}_{QNR}$).

Example: for $p = 3$, we have $\mathcal{G}_S = \mathcal{G}_{QR} = \{1\}$ and $\mathcal{G}_Z = \mathcal{G}_{QNR} = \{2\}$.

**Property 1.2**: Let $p \in \mathbb{P} \setminus \{2\}$. We have $\mathcal{G}_Z \subset \mathcal{G}_{QNR}$ and $\mathcal{G}_S \subset \mathcal{G}_{QR}$.

Proof: The second inclusion is immediate, since a quadratic residue is defined as a square modulo $p$, and a semi-primitive root is the square of a primitive root. Moreover, $g \in \mathbb{F}_p^*$ is either a quadratic non-residue or a quadratic residue, i.e. if $h$ is not a quadratic non-residue, there exists $k$ such that $g = h^{2k}$. Thus, the order of $\langle g \rangle$ divides $\frac{p-1}{2}$ and $g$ is therefore not a primitive root, hence $\mathcal{G}_Z \subset \mathcal{G}_{QNR}$.

**Definition 1.3**: Let $p \in \mathbb{P} \setminus \{2\}$ and $m \in \mathbb{Z}/p\mathbb{Z}$. We denote by $\mathcal{R}_p(m)$ the set of square roots of $m$ i.e.:

$$\mathcal{R}_p(m) = \{g \in \mathbb{Z}/p\mathbb{Z} | g^2 \equiv m\}.$$

More generally, we recursively define the sequence $\left(\mathcal{R}_p^k(m)\right)_{k \geq 0}$ by:

$$\mathcal{R}_p^0(m) = \{m\}, \qquad \mathcal{R}_p^{k+1}(m) = \bigcup_{g \in \mathcal{R}_p^k(m)} \mathcal{R}_p(g)$$

so that $\mathcal{R}_p(m) = \mathcal{R}_p^1(m)$.

We will simply note $\mathcal{R}(m)$ and $\mathcal{R}^k(m)$ when the value of $p$ is unambiguous.

We will also note $\mathcal{R}^k$ for $\mathcal{R}^k(1)$.

**Property 1.3**: $\mathcal{R}^k$ is a subgroup of $\mathbb{F}_p^*$.

**Property 1.4**: let $p \in \mathbb{P} \setminus \{2\}$ and $n$ the highest power of 2 dividing $p - 1$. For any $0 \leq k \leq n$, we have:

$$|\mathcal{R}^k| = 2^k$$

and for any $1 \leq k \leq n$:

$$\mathcal{R}^{k-1} \subset \mathcal{R}^k, \qquad |\mathcal{R}^k \setminus \mathcal{R}^{k-1}| = 2^{k-1}$$

For $k > n$, we have $\mathcal{R}^k = \mathcal{R}^n$.

Proof: If $h$ is a primitive root, we note that $\mathcal{R}^k = \left\{h^k, k \in \frac{p-1}{2^k}\mathbb{Z}\right\}$. The first two points follow from this. Furthermore, for $k > n$, $\mathbb{F}_p^*$ contains no element of order $2^k$, thus $\mathcal{R}^k = \mathcal{R}^n$.

## 2. Generators and semi-generators of $\mathbb{F}_p^*$

In this section, we show that $\mathcal{G}_Z$ can be obtained from quadratic residues. We give an algorithm that generates $\mathcal{G}_Z$ without computing any $gcd$.

## 2.1. Determining primitives

Let $p \in \mathbb{P} \setminus \{2\}$. The number of generators of the cyclic group $\mathbb{F}_p^*$ is equal to $\varphi(\varphi(p)) = \varphi(p-1)$. The following property gives a constructive method for determining primitive or semi-primitive roots.

**Proposition 2.1**: We decompose into prime factors $p - 1 = \prod_{q \in \mathbb{P}} q^{\alpha_q}$. Then $g \in \mathbb{F}_p^*$ is a primitive root if and only if:

$$\forall q \in \mathbb{P} \ \left(\alpha_q \geq 1 \Rightarrow g^{\frac{p-1}{q}} \not\equiv 1\right)$$

In this case we have:
$$\mathcal{G}_Z = \{g^k, 1 \leq k \leq p-2, \gcd(k, p-1) = 1\}$$

Moreover if $p \in \mathcal{P}_{1,4}$:

$$g \in \mathcal{G}_S \Leftrightarrow g \in \mathcal{G}_{QR} \text{ and } \forall q \in \mathbb{P} \ \left(\alpha_q \geq 1 \Rightarrow g^{\frac{p-1}{2q}} \not\equiv 1\right) \Leftrightarrow \emptyset \subsetneq \mathcal{R}_p(g) \subset \mathcal{G}_Z.$$

If $p \in \mathcal{P}_{3,4}$:

$$g \in \mathcal{G}_S \Leftrightarrow g^{\frac{p-1}{2}} \equiv 1 \text{ and } \forall q \in \mathbb{P} \setminus \{2\} \ \left(\alpha_q \geq 1 \Rightarrow g^{\frac{p-1}{2q}} \not\equiv 1\right)$$

In this second case, only one of the two square roots of any $g \in \mathcal{G}_S$ is a primitive root.

**Proof**: Let us take $g \in \mathcal{G}_Z$. We know that the order of $g^k$ is $\frac{p-1}{\gcd(k, p-1)}$. More specifically, for any $q \in \mathbb{P}$ such that $\alpha_q \geq 1$ (i.e. $q$ dividing $p - 1$) the order of $g^{\frac{p-1}{q}}$ is $\frac{p-1}{\gcd\left(\frac{p-1}{q}, p-1\right)} = q$ so $g^{\frac{p-1}{q}} \not\equiv 1$. Furthermore, $h \in \mathcal{G}_Z$ if and only if there is $k$ coprime with $p - 1$ such that $h = g^k$.

Conversely, let $g \in \mathbb{F}_p^*$ such that for any prime $q$ dividing $p - 1$, $g^{\frac{p-1}{q}} \not\equiv 1$. The order of $g$ divides $p - 1$. However, the hypothesis implies that it cannot be a strict divisor of $p - 1$, thus $g$ is a primitive root.

We now assume that $p \in \mathcal{P}_{1,4}$. By definition, $g \in \mathcal{G}_S \Leftrightarrow \exists h \in \mathcal{G}_Z \ g = h^2$, i.e. $g$ is in $\mathcal{G}_S$ if and only if there exists $h$ such that $g = h^2$ and $\forall q \in \mathbb{P} \ \left(\alpha_q \geq 1 \Rightarrow h^{\frac{p-1}{q}} \not\equiv 1\right)$. Since $\frac{p-1}{q}$ is always even (including when $q = 2$ because $p - 1$ is a multiple of 4), we have $h^{\frac{p-1}{q}} = g^{\frac{p-1}{2q}}$ hence:

$$g \in \mathcal{G}_S \Leftrightarrow g \in \mathcal{G}_{QR} \text{ and } \forall q \in \mathbb{P} \ \left(\alpha_q \geq 1 \Rightarrow g^{\frac{p-1}{2q}} \not\equiv 1\right).$$

Finally, $g \in \mathcal{G}_S$ if and only it has two square roots $h$ and $-h$, one of which at least (say $h$) is an element of $\mathcal{G}_Z$. But then as for any prime $q$ dividing $p - 1$, $\frac{p-1}{q}$ is even, we have $(-h)^{\frac{p-1}{q}} = h^{\frac{p-1}{q}} \not\equiv 1$, hence $-h \in \mathcal{G}_Z$ too.

Finally, when $p \in \mathcal{P}_{3,4}$, $\frac{p-1}{2}$ is odd and $-1$ is not a quadratic residue. If $g \in \mathcal{G}_S$ then as in the previous paragraph, $g \in \mathcal{G}_{QR}$ and $\forall q \in \mathbb{P} \setminus \{2\}$ $\left(\alpha_q \geq 1 \Rightarrow g^{\frac{p-1}{2q}} \not\equiv 1\right)$. Conversely, assume $g = h^2$ and $\forall q \in \mathbb{P} \setminus \{2\}$ $\left(\alpha_q \geq 1 \Rightarrow g^{\frac{p-1}{2q}} \not\equiv 1\right)$. Then $\forall q \in \mathbb{P} \setminus \{2\}$ $\left(\alpha_q \geq 1 \Rightarrow (\pm h)^{\frac{p-1}{q}} \not\equiv 1\right)$ and $h^{\frac{p-1}{2}} \in \{\pm 1\}$ thus $h^{\frac{p-1}{2}} = -1$ or $(-h)^{\frac{p-1}{2}} = -1$. In all cases, either $h$ or $-h$ is a primitive root, but never both.

## 2.2. Relationships between primitive and semi-primitive roots

### a. Primitive roots obtained from residues

**Definition 2.2.1**: if $h$ is an element of order $d$ in a group and $0 \leq x \leq d$ we note $\langle h \rangle^{(x)} = \{h, \ldots, h^x\}$. It is a subset of the subgroup generated by $h$ of cardinal $x$.

In this whole section, we let $p \in \mathbb{P} \setminus \{2\}$ and we write $p - 1 = 2^n z$ with $z$ odd.

**Proposition 2.2.1**: Let $g \in \mathcal{G}_{QNR}$. For any $0 \leq t \leq n$, the multiplicative order of $g^{2^t z}$ is $2^{n-t}$. In other words, $\langle g^{2^t z} \rangle = \mathcal{R}^{(n-t)}$. Moreover, if we note $b = g^{2^t z}$, with $0 < t \leq n$, for any integer $x$, $b^{2x} \in \mathcal{R}^{(n-t-1)}$ and $b^{2x+1} \in \mathcal{R}^{(n-t)} \setminus \mathcal{R}^{(n-t-1)}$.

In particular, when $t = 0$, the odd powers of $b$ are in $\mathcal{G}_{QNR}$ and the even powers in $\mathcal{G}_{QR}$.

*Proof*: The order of $g^z$ divides $2^n$. By Euler's criterion, $(g^z)^{\frac{p-1}{2}} \equiv (-1)^z = -1$, thus $g^z$ is in $\mathcal{G}_{QNR}$ and a generator of $\mathcal{R}^{(n)}$. Hence, $g^{2^t z}$ is a generator of $\mathcal{R}^{(n-t)}$. The other assertions are immediate.

**Corollary 2.2.1**: Let $m$ be a semi-primitive root. Then:
$$\mathcal{G}_{QR} = \langle m \rangle = \langle m \rangle^{(z)} \mathcal{R}^{(n-1)}.$$

More generally, generators of $\mathcal{G}_{QR}$ coincide with semi-primitive roots.

*Proof*: We know that $\mathcal{G}_{QR}$ is the subgroup made of the roots of the polynomial $\left(X^{\frac{p-1}{2}} - 1\right)$. By definition, a semi-primitive root is of order equal to $\frac{p-1}{2}$ and thus a generator of this subgroup, and conversely. Furthermore, we have $\langle m \rangle^{(z)} \subset \mathcal{G}_{QR}$ and $\mathcal{R}^{(n-1)} \subset \mathcal{G}_{QR}$ hence $\langle m \rangle^{(z)} \mathcal{R}^{(n-1)} \subset \mathcal{G}_{QR}$.

Moreover, $\left|\langle m \rangle^{(z)}\right| = z$, $\left|\mathcal{R}^{(n-1)}\right| = 2^{n-1}$. If $m^k r \equiv m^{k'} r'$ with $1 \leq k \leq k' \leq z$ and $r, r' \in \mathcal{R}^{(n-1)}$ thus $m^{k'-k} r' \equiv r$ hence $m^{2^{n-1}(k'-k)} \equiv 1$. This implies $k = k'$ and $r \equiv r'$ follows. Thus, we have $\left|\langle m \rangle^{(z)}\right| \cdot \left|\mathcal{R}^{(n-1)}\right| = 2^{n-1} z = \frac{p-1}{2}$ distinct residues in $\langle m \rangle^{(z)} \mathcal{R}^{(n-1)}$ hence $\langle m \rangle^{(z)} \mathcal{R}^{(n-1)} = \mathcal{G}_{QR}$.

**Proposition 2.2.2**: For $m$ a semi-primitive root, we have:
$$m\left(\mathcal{R}^{(n)} \setminus \mathcal{R}^{(n-1)}\right) \subset \mathcal{G}_Z.$$

*Proof*: Let $r \in \mathcal{R}^{(n)} \setminus \mathcal{R}^{(n-1)}$, and let us show that $mr \in \mathcal{G}_Z$. Let $q$ be a prime factor of $p - 1 = 2^n z$.

If $q = 2$, we have $(mr)^{\frac{p-1}{q}} \equiv m^{2^{n-1}z} r^{2^{n-1}z} \equiv r^{2^{n-1}z} \not\equiv 1$ because $r^z \in \mathcal{R}^{(n)} \setminus \mathcal{R}^{(n-1)}$ since $z$ is odd.

If $q \neq 2$, we have $(mr)^{\frac{p-1}{q}} \equiv m^{\frac{2^n z}{q}} r^{\frac{2^n z}{q}} \equiv m^{\frac{2^n z}{q}} \not\equiv 1$ because $2^{n-1}z$ does not divide $\frac{2^n z}{q}$.

Thus, using proposition 2.1, we deduce that $mr \in \mathcal{G}_Z$.

**Proposition 2.2.3**: If $p \in \mathcal{P}_{3,4}$, there are as many semi-primitive roots as primitive roots, and if $p \in \mathcal{P}_{1,4}$ there are half as many.

*Proof*: We know that there are $\varphi(p-1)$ primitive roots and, by corollary 2.2.1, $\varphi\left(\frac{p-1}{2}\right)$ semi-primitive roots modulo $p$. We also know that $\varphi(ab) = \varphi(a)\varphi(b)$ if $a, b$ are coprime and $\varphi(2^k) = 2^{k-1}$.

If $p \equiv 3[4]$ then $n = 1$ i.e. $\frac{p-1}{2}$ is odd and $\varphi(p-1) = \varphi(2)\varphi\left(\frac{p-1}{2}\right) = \varphi\left(\frac{p-1}{2}\right)$.

If $p \equiv 1[4]$ then $n \geq 2$ and $\varphi(p-1) = \varphi(2^n)\varphi(z) = 2^{n-1}\varphi(z)$ whereas $\varphi\left(\frac{p-1}{2}\right) = \varphi(2^{n-1})\varphi(z) = 2^{n-2}\varphi(z)$.

**Proposition 2.2.4**: We assume $p \in \mathcal{P}_{1,4}$. Let $m$ be a semi-primitive root. Then:
$$-m\big(\mathcal{R}^{(n)} \setminus \mathcal{R}^{(n-1)}\big) \subset \mathcal{G}_Z$$
Moreover, $-m \in \mathcal{G}_S$ if and only if $p \in \mathcal{P}_{1,8}$.

*Proof*: Let $r \in \mathcal{R}^{(n)} \setminus \mathcal{R}^{(n-1)}$, and let us show that $-mr \in \mathcal{G}_Z$ like in proposition 2.2.2. Let $q$ be a prime factor of $p - 1 = 2^n z$. As $p \in \mathcal{P}_{1,4}$, we have $n \geq 2$ so in all cases $(-mr)^{\frac{p-1}{q}} \equiv (mr)^{\frac{p-1}{q}} \not\equiv 1$, because $mr \in \mathcal{G}_Z$. Therefore $-mr \in \mathcal{G}_Z$.

Moreover $-m \in \mathcal{G}_S$ if and only if $(-m)^{\frac{p-1}{2q}} \not\equiv 1$ for any prime factor $q$ of $p - 1$, which is true when $n \geq 3$, because then $(-m)^{\frac{p-1}{2q}} \equiv m^{\frac{p-1}{2q}}$. Otherwise $m^{\frac{p-1}{4}} \equiv -1$ thus $(-m)^{\frac{p-1}{4}} \equiv 1$ i.e. $-m \notin \mathcal{G}_S$. We hence proved that $-m \in \mathcal{G}_S$ if and only if $p \in \mathcal{P}_{1,8}$.

**Proposition 2.2.5**: Let $m$ be a semi-primitive root. We first assume $p \in \mathcal{P}_{3,4}$. Then $m^2 \in \mathcal{G}_S$ and in particular:
$$m^2\big(\mathcal{R}^{(n)} \setminus \mathcal{R}^{(n-1)}\big) \subset \mathcal{G}_Z$$
We now assume $p \in \mathcal{P}_{1,4}$. Then:
$$\pm m^2\big(\mathcal{R}^{(n)} \setminus \mathcal{R}^{(n-1)}\big) \subset \mathcal{G}_Z$$
if $p \in \mathcal{P}_{1,8}$ then $\pm m^2 \notin \mathcal{G}_S$. If $p \in \mathcal{P}_{5,8}$ then $m^2 \notin \mathcal{G}_S$ and $-m^2 \in \mathcal{G}_S$. In the last case, we identify $-\mathcal{G}_S$ to the set of squares of $\mathcal{G}_S$.

*Proof*: When $p \in \mathcal{P}_{3,4}$, $n = 1$ so $\mathcal{R}^{(n)} \setminus \mathcal{R}^{(n-1)} = \{-1\}$. By proposition 2.2.1 we have $-m \in \mathcal{G}_Z$ hence $m^2 \in \mathcal{G}_S$.

When $p \in \mathcal{P}_{1,4}$ i.e. $n \geq 2$, let us take $r \in \mathcal{R}^{(n)} \setminus \mathcal{R}^{(n-1)}$ and focus on $(\pm m^2 r)^{\frac{p-1}{q}} = (m^2 r)^{\frac{p-1}{q}}$ with $q$ dividing $p-1$.

If $q = 2$ we have $(m^2 r)^{\frac{p-1}{2}} = r^{2^{n-1}z} \not\equiv 1$ because $r^z \in \mathcal{R}^{(n)} \setminus \mathcal{R}^{(n-1)}$.

If $q \neq 2$ we have $(m^2 r)^{\frac{p-1}{q}} = m^{\frac{2(p-1)}{q}} \not\equiv 1$ because $\frac{p-1}{2}$ does not divide $\frac{2(p-1)}{q}$.

We thus proved that $\pm m^2 (\mathcal{R}^{(n)} \setminus \mathcal{R}^{(n-1)}) \subset \mathcal{G}_Z$.

If $p \in \mathcal{P}_{1,8}$ then $(\pm m^2)^{\frac{p-1}{4}} \equiv 1$ and $\pm m^2 \notin \mathcal{G}_S$.

If $p \in \mathcal{P}_{5,8}$ then $n = 2$ and $(m^2)^{\frac{p-1}{4}} \equiv 1$, $(-m^2)^{\frac{p-1}{4}} \equiv (-1)^z (m^2)^{\frac{p-1}{4}} \equiv -1$ thus $m^2 \notin \mathcal{G}_S$ and $-m^2 \in \mathcal{G}_S$ (checking the other powers is immediate). Similarly, if $m \in \mathcal{G}_S$, $-m \notin \mathcal{G}_S$ thus the square function is injective on $\mathcal{G}_S$, which shows that $-\mathcal{G}_S$ is equal to the set of squares of $\mathcal{G}_S$.

Propositions 2.2.4 and 2.2.5 allow to build primitive roots from semi-primitive roots and elements of $\mathcal{R}^{(n)} \setminus \mathcal{R}^{(n-1)}$.

**Definition 2.2.2**: For some $r \in \mathcal{R}^{(n)} \setminus \mathcal{R}^{(n-1)}$, we define the set $\mathcal{G}'_Z \coloneqq r\mathcal{G}_Z$.

**Proposition 2.2.6**: $\mathcal{G}'_Z$ is independent of the choice of $r$ and contains $\mathcal{G}_S$. It has <u>same size</u> as $\mathcal{G}_Z$ and we have $\mathcal{G}_Z = r\mathcal{G}'_Z$ for any $r \in \mathcal{R}^{(n)} \setminus \mathcal{R}^{(n-1)}$. Moreover, if for $k \in \mathbb{N}^*$ we define $\mathcal{G}_S^{(k)} = \{m^k, m \in \mathcal{G}_S\}$:

- when $p \in \mathcal{P}_{3,4}$, $\mathcal{G}'_Z = \mathcal{G}_S = \mathcal{G}_S^{(2)}$ whereas $\mathcal{G}_Z = -\mathcal{G}'_Z$.
- when $p \in \mathcal{P}_{5,8}$, $\mathcal{G}'_Z = \mathcal{G}_S \sqcup -\mathcal{G}_S = \mathcal{G}_S \sqcup \mathcal{G}_S^{(2)}$.
- when $p \in \mathcal{P}_{9,16}$, $\mathcal{G}'_Z = \mathcal{G}_S \sqcup \left(\mathcal{G}_S^{(2)} \sqcup -\mathcal{G}_S^{(2)}\right)$
- for any $p$, we have:

$$\mathcal{G}'_Z = \bigcup_{t=0}^{n-1} \mathcal{G}_S^{(2^t)} = \left\{g \in \mathbb{F}_p^* \middle| \forall q \in \mathbb{P} \setminus \{2\} \left(q | p-1 \Rightarrow g^{\frac{p-1}{2q}} \not\equiv 1\right)\right\}.$$

*Proof*: To prove that $\mathcal{G}'_Z$ does not depend on $r$, it is enough to observe that if $r, r' \in \mathcal{R}^{(n)} \setminus \mathcal{R}^{(n-1)}$ and $g \in \mathcal{G}_Z$, we have $rr'g \in \mathcal{G}_Z$. This can be done using Proposition 2.1: we know that $r^{\frac{p-1}{2}} \equiv r'^{\frac{p-1}{2}} \equiv -1$ whereas if $q$ is an odd prime factor of $p-1$, $r^{\frac{p-1}{q}} \equiv r'^{\frac{p-1}{q}} \equiv 1$. Proposition 2.2.2 ensures that $\mathcal{G}_S \subset \mathcal{G}'_Z$. It is also clear that $\mathcal{G}'_Z$ and $\mathcal{G}_Z$ have the same size and that $\mathcal{G}_Z = r\mathcal{G}'_Z$ for any $r \in \mathcal{R}^{(n)} \setminus \mathcal{R}^{(n-1)}$.

Finally:

- If $p \in \mathcal{P}_{3,4}$, since from proposition 2.2.6 the two sets have the same size, $\mathcal{G}'_Z = \mathcal{G}_S$. Proposition 2.2.5 also ensures that $\mathcal{G}_S^{(2)} \subset \mathcal{G}_S$ and since $-1 \notin \mathcal{G}_S$ the square function is injective from $\mathcal{G}_S$ to $\mathcal{G}_S^{(2)}$ thus $\mathcal{G}_S^{(2)} = \mathcal{G}_S$. In this case we know that $\mathcal{R}^{(n)} \setminus \mathcal{R}^{(n-1)} = \{-1\}$, proposition 2.2.2 ensures $\mathcal{G}_Z = -\mathcal{G}'_Z$.

- If $p \in \mathcal{P}_{5,8}$, proposition 2.2.5 ensures that $-\mathcal{G}_S = \mathcal{G}_S^{(2)}$ is disjoint of $\mathcal{G}_S$ and included in $\mathcal{G}'_Z$ hence a cardinality argument ensures $\mathcal{G}'_Z = \mathcal{G}_S \sqcup -\mathcal{G}_S = \mathcal{G}_S \sqcup \mathcal{G}_S^{(2)}$.
- If $p \in \mathcal{P}_{1,8}$, proposition 2.1 shows that $\mathcal{G}_S = -\mathcal{G}_S$ thus the size of $\mathcal{G}_S^{(2)}$ is half that of $\mathcal{G}_S$. Proposition 2.2.5 ensures $\mathcal{G}_S \sqcup \left(\mathcal{G}_S^{(2)} \cup -\mathcal{G}_S^{(2)}\right) \subset \mathcal{G}'_Z$. Remains to show that $\mathcal{G}_S^{(2)}$ is disjoint of $-\mathcal{G}_S^{(2)}$ if and only if $p \in \mathcal{P}_{9,16}$ i.e. $n = 3$. If $\mathcal{G}_S^{(2)}$ is not disjoint of $-\mathcal{G}_S^{(2)}$, there exists $g \in \mathcal{G}_Z$ and $k$ coprime to $p-1$ such that $g^{4k} \equiv -g^4$, hence, by taking odd powers coprime to $p-1$, $\mathcal{G}_S^{(2)} = -\mathcal{G}_S^{(2)}$ while on the other hand $\left(g^{\frac{k-1}{2}}\right)^8 \equiv -1 \equiv (g^z)^{2^{n-1}}$ which implies $n \geq 4$. Conversely if $n \geq 4$, then $-1$ has a square root which leaves $\mathcal{G}_S$ stable, thus $\mathcal{G}_S^{(2)} = -\mathcal{G}_S^{(2)}$.
- In general, similarly to previous items, for any $t \leq n-1$, $\mathcal{G}_S^{(2^t)}$ is the set of elements of order $\frac{p-1}{2^{t+1}}$, and if moreover $t \neq n-1$ then $\mathcal{G}_S^{(2^t)} = -\mathcal{G}_S^{(2^t)}$ thus by induction $\left|\mathcal{G}_S^{(2^t)}\right| = \frac{1}{2^t}|\mathcal{G}_S| = \frac{1}{2^{t+1}}|\mathcal{G}_Z|$. For $t = n-1$ we have $-\mathcal{G}_S^{(2^{n-1})} \cap \mathcal{G}_S^{(2^{n-1})} = \emptyset$ hence $\left|\mathcal{G}_S^{(2^{n-1})}\right| = \frac{1}{2^{n-1}}|\mathcal{G}_Z|$. We deduce that the sets $\left(\mathcal{G}_S^{(2^t)}\right)_{0 \leq t \leq n-1}$ are pairwise disjoint and their union has the same cardinal as $\mathcal{G}_Z$. Proposition 2.1 shows that if $r \in \mathcal{R}^{(n)} \setminus \mathcal{R}^{(n-1)}$, $r\mathcal{G}_S^{(2^t)} \subset \mathcal{G}_Z$. We must then have $\mathcal{G}'_Z = \cup_{t=0}^{n-1} \mathcal{G}_S^{(2^t)}$. Proposition 2.1 also shows that $\mathcal{G}'_Z = \left\{g \in \mathbb{F}_p^* \,\middle|\, \forall q \in \mathbb{P} \setminus \{2\} \left(q|p-1 \Rightarrow g^{\frac{p-1}{2q}} \not\equiv 1\right)\right\}$.

**Theorem 1**: We have:

- $g \in \mathcal{G}_Z$ if and only if there exists $m \in \mathcal{G}_{QR}$ such that $\forall q \in \mathbb{P} \setminus \{2\}\, q|p-1 \Rightarrow m^{\frac{p-1}{2q}} \not\equiv 1$ and $r \in \mathcal{R}^{(n)} \setminus \mathcal{R}^{(n-1)}$ such that $g = mr$.
- $g \in \mathcal{G}_Z$ if and only if $g \in \mathcal{G}_{QNR}$ and $\forall q \in \mathbb{P} \setminus \{2\}\, q|p-1 \Rightarrow g^{\frac{p-1}{q}} \not\equiv 1$.

*Proof*: We just need to prove the first item, the second one has already been established in proposition 2.1.

The Chinese theorem identifies $\mathbb{F}_p^*$ to $(\mathbb{Z}/2^n\mathbb{Z}) \times (\mathbb{Z}/z\mathbb{Z})$. A primitive root $g$ corresponds to a pair $(x, y)$ with $x$ odd and $y$ coprime to $z$. Let $m$ be the element corresponding to $(2x, y)$ and $r$ the element corresponding to $(-x, 0)$. Then it is clear that $g = mr$ but $r$ is of order $2^n$ so in $\left(\mathcal{R}^{(n)} \setminus \mathcal{R}^{(n-1)}\right)$, while the order of $m$ is $2^{n-1}z$ i.e. $m$ is a semi-primitive root and in particular $m \in \mathcal{G}_{QR}$ hence $\forall q \in \mathbb{P} \setminus \{2\}\, q|p-1 \Rightarrow m^{\frac{p-1}{2q}} \not\equiv 1$.

Conversely, using this identification, $m$ corresponds to $(x, y)$ with $x$ even and $y$ coprime to $z$, and $r$ corresponds to $(t, 0)$ with $t$ odd. Thus $mr$ corresponds to $(x + t, y)$ which is indeed a generator of $(\mathbb{Z}/2^n\mathbb{Z}) \times (\mathbb{Z}/z\mathbb{Z})$ because $x + t$ is odd and y coprime to $z$.

### b. A sequence of primitive and semi-primitive roots

Knowing a single primitive root allows to obtain every element of $\mathcal{G}_Z$ by exponentiation. Knowing a semi-primitive root and an element of $\mathcal{R}^{(n)} \setminus \mathcal{R}^{(n-1)}$ also allows to find one thus every primitive root. We describe here an alternative method for generating primitive roots.

**Corollary 2.2.5**: Let $m \in \mathcal{G}_S$ and $r \in \mathcal{R}^{(n)} \setminus \mathcal{R}^{(n-1)}$ fixed, then:

$$\forall k \in \mathbb{N}^* \ \forall a \in \mathbb{P} \setminus \mathcal{D}_p(z) \ m^{a^k}r \in \mathcal{G}_z.$$

*Proof*: We apply the first item of theorem 1. If $q \in \mathbb{P} \setminus \{2\}$ divides $p-1$, then $m^{\frac{a^k(p-1)}{2q}} \not\equiv 1$ otherwise $\frac{p-1}{2}$ would divide $\frac{a^k(p-1)}{2q}$, hence by Gauss's theorem $\frac{p-1}{2}$ would divide $\frac{p-1}{2q}$ which is impossible.

**Definition 2.2.3**: Let $a \in \mathbb{P}$ coprime to $z$, the underline{discrete logarithm of 1 to the base $a$ modulo $z$}, denoted by $\text{logd}_a(1,z)$, is the multiplicative order of $a$ modulo $z$, i.e. the smallest integer $k$ such that $a^k \equiv 1 \ [z]$.

**Proposition 2.2.7**: Let $m \in \mathcal{G}_S$, $b = m^{2^{n-1}}$ and $k = \text{logd}_2(1,z)$. Then:

$$b^{2^k} \equiv b[p].$$

*Proof*: By definition $2^k \equiv 1 \ [z]$. We also know that $b = m^{2^{n-1}}$ is of order $z$ modulo $p$, so $b^{2^k} \equiv b^1 \equiv b$.

**Remark 2.2.1**:

- There are several algorithms to calculate the discrete logarithm, all of which take sub-exponential time. See for example Pollard's rho algorithm, Pohlig-Hellman algorithm or the general number field sieve (GNFS).
- Searching for the discrete logarithm to the base 2 is adapted to the binary system used by the computer

**Proposition 2.2.8**: Let $a$ be coprime to $z$. Then we have:

$$\text{logd}_a(1,z) \leq \operatorname*{lcm}_{q \in \mathcal{D}_p(p-1)} q^{\alpha_q - 1}(q-1).$$

*Proof*: The Chinese theorem gives $\mathbb{Z}/z\mathbb{Z} \simeq \prod_{q \in \mathcal{D}_p(p-1)} \mathbb{Z}/q^{\alpha_q}\mathbb{Z}$. Therefore the unit group $(\mathbb{Z}/z\mathbb{Z})^\times$ can be identified to $\prod_{q \in \mathcal{D}_p(p-1)}(\mathbb{Z}/q^{\alpha_q}\mathbb{Z})^\times$ and any element of this group is of order at most $\operatorname*{lcm}_{q \in \mathcal{D}_p(p-1)}|(\mathbb{Z}/q^{\alpha_q}\mathbb{Z})^\times| = \operatorname*{lcm}_{q \in \mathcal{D}_p(p-1)} q^{\alpha_q-1}(q-1)$. The same applies to $\text{logd}_a(1,z)$ which is the multiplicative order of $a$ in $\mathbb{Z}/z\mathbb{Z}$.

**Proposition 2.2.9**: If $m \in \mathcal{G}_Z$ (respectively $\mathcal{G}_S, \mathcal{G}_{QR}, \mathcal{G}_{NQR}$) then $m^{p-2} \in \mathcal{G}_Z$ (respectively $\mathcal{G}_S$, $\mathcal{G}_{QR}, \mathcal{G}_{NQR}$).

*Proof*: $m^{p-2}$ is the multiplicative inverse of $m$ in $\mathbb{F}_p^*$. As the inverse of a generator is a generator, $\mathcal{G}_Z$ and $\mathcal{G}_S$ are stable by conversion to the inverse. Moreover, $\mathcal{G}_{QR}$ is a multiplicative subgroup and $\mathcal{G}_{NQR}$ its complementary, so it's the same for these two sets.

**Corollary 2.2.9**: For $t \in \{0, \ldots, n-1\}$, we let $a_t = 2^{n-t}z - 1$. For any $m \in \mathcal{G}_Z$, $m^{2a_0} \equiv m^{2a_1}, m^{2a_2}, \ldots, m^{2a_{n-1}}$ are $n-1$ distinct semi-primitives roots.

*Proof*: Clearly $a_t$ is odd and coprime to $z$, thus also to $p-1$. We deduce that $m^{a_t} \in \mathcal{G}_Z$ hence $m^{2a_t} \in \mathcal{G}_S$. Moreover, $2 \leq 2 \times (2z-1) = 2a_{n-1} < \cdots < 2a_1 = 2 \times (2^{n-1}z - 1) = p - 3$ so $m^{2a_1} \ldots m^{2a_{n-1}}$ are pairwise distinct. Finally, $2a_0 = 2p - 4 = p - 1 + p - 3$ and $m^{p-1} \equiv 1$ therefore we have $m^{2a_0} \equiv m^{2a_1}$.

The following proposition gives a family of sequences of semi-generators and generators, as well as a method to obtain each element of the sets $\mathcal{G}_S$ and $\mathcal{G}_Z$.

**Proposition 2.2.10**: Let $g \in \mathcal{G}_S$ and $m' \in \mathcal{R}^{(n)} \setminus \mathcal{R}^{(n-1)}$ fixed such that $g^z m'^{2z} \equiv 1$. We define the sequence $(U_x)$ by:

$$U_0 = g, \qquad U_{x+1} = (m'U_x)^2$$

For all $x \in \mathbb{N}$, $U_x \in \mathcal{G}_S$ and $m'U_x \in \mathcal{G}_Z$. Moreover if $0 \leq x < y < \text{logd}_2(1,z)$, $U_x \not\equiv U_y$ and $(U_x)$ is periodic of period $\text{logd}_2(1,z)$.

*Proof*: In the identification of $\mathbb{F}_p^*$ to $(\mathbb{Z}/2^n\mathbb{Z}) \times (\mathbb{Z}/z\mathbb{Z})$, any element of $\mathcal{G}_S$ corresponds to a pair $(2k, l)$ with $k$ odd and $l$ coprime to $z$ while any elements of $\mathcal{R}^{(n)} \setminus \mathcal{R}^{(n-1)}$ correspond to a pair $(i, 0)$ with $i$ odd. From this and the equation $g^z m'^{2z} \equiv 1$, we deduce that there is $k_0$ odd and $l_0$ coprime to $z$ such that $g$ corresponds to $(2k_0, l_0)$ and $m'$ to $(-k_0, 0)$.

Thus, by recurrence, we note that $U_x$ corresponds to $(2k_x, l_x) = (2k_0, 2^x l_0)$ where $2^x l_0$ is still coprime to $z$, so $U_x \in \mathcal{G}_S$. Similarly, $m'U_x$ corresponds to $(k_0, 2^x l_0)$ and we have $m'U_x \in \mathcal{G}_Z$. Moreover, by definition of $\text{logd}_2(1,z)$, the elements of $U_x$ are pairwise distinct for $0 \leq x < \text{logd}_2(1,z)$ and $U_{\text{logd}_2(1,z)} = U_0$ which implies the periodicity of $(U_x)$.

**Remark 2.2.2**: We verify that $U_x = g^{2^x} m'^{2(2^x-1)}$. In general, the elements $(m'U_x)$ form a strict subset of $\mathcal{G}_Z$. However, for any generator $h \in \mathcal{G}_Z$ there exists $a$ coprime to $p-1$ such that $h = (m'g)^a$. We can also choose $a$ to be prime. Thus, if we take $a \in \mathbb{P} \setminus \mathcal{D}(p-1)$ such that $g^a m'^a$ has not already been generated then the $\text{logd}_2(1,z)$ first terms $U_x^{(a)} := g^{2^x a} m'^{2a(2^x-1)}$ will all be new elements of $\mathcal{G}_Z$. We deduce that there exists a set $A$ of $\frac{\varphi(p-1)}{\text{logd}_2(1,z)}$ primes, odd and coprime to $z$, and such that:

$$\mathcal{G}_Z = \left\{ g^{2^x a} m'^{(2^{x+1}-1)a}, (a, x) \in (\{1\} \cup A) \times [\![0, \text{logd}_2(1,z) - 1]\!] \right\}$$

*Proof*: By construction we have $U_x = g^{2^x} m'^{2 + \cdots + 2^x} = g^{2^x} m'^{2(2^x-1)}$. If $h \in \mathcal{G}_Z$, we know that $m'g \in \mathcal{G}_Z$ so there exists $a$ coprime to $p-1$ (i.e. odd and coprime to $z$) such that $h = (m'g)^a$. Dirichlet's theorem ensures that we can choose $a$ to be prime, even if it means adding it a multiple of $p-1$.

Identifying $\mathbb{F}_p^*$ to $(\mathbb{Z}/2^n\mathbb{Z}) \times (\mathbb{Z}/z\mathbb{Z})$ as in the proof of the previous results, $m'^a U_x^{(a)}$ correspond to $(ak_0, 2^x al_0)$. Let $a, b, x, y$ be such that $(ak_0, 2^x al_0) = (bk_0, 2^y bl_0)$. In particular $a.(-k_0, 0) = b.(-k_0, 0)$ thus $m'^a \equiv m'^b$ hence for any $k \in \mathbb{N}$, $U_{x+k}^{(a)} = U_{y+k}^{(b)}$ from which it follows by periodicity of these sequences that $m'^b U_0^{(b)}$ is equal to one of the terms of $\left( m'^a U_x^{(a)} \right)$. Therefore, if we construct $A$ recursively in such a way that $m'^a U_0^{(a)}$ is not among the generators already constructed, we keep adding chunks of $\text{logd}_2(1,z)$ generators, which will eventually cover the entire set $\mathcal{G}_Z$.

**Remark 2.2.3**: For each value of $a \in \{1\} \cup A$, we obtain $\text{logd}_2(1,z)$ pairwise distinct generators. We deduce that $\text{logd}_a(1,z)$ divides $|\mathcal{G}_Z| = \varphi(p-1) = \varphi(2^n)\varphi(z)$. In fact, since $\text{logd}_2(1,z)$ is the (multiplicative) order of 2 modulo $z$, we know that it divides $\varphi(z)$, the order of the group $(\mathbb{Z}/z\mathbb{Z})^\times$. Remark also that $\logd_2(1,z) > \log_2(z)$ since necessarily $2^{\logd_2(1,z)} > z$.

**Corollary 2.2.4**: We have $\logd_2(1,z) \leq \varphi(z) \leq \left\lfloor \frac{z}{\logd_2(1,z)} \right\rfloor \logd_2(1,z)$.

*Proof*: As stated in the remark 2.2.3, we know that $\logd_2(1,z)$ divides $\varphi(z)$. Therefore $\logd_2(1,z) \leq \varphi(z)$ and $\varphi(z) = \frac{\varphi(z)}{\logd_2(1,z)} \logd_2(1,z) \leq \left\lfloor \frac{z}{\logd_2(1,z)} \right\rfloor \logd_2(1,z)$.

**Proposition 2.2.11**: Let $g$ and $m'$ be as in proposition 2.2.10. There exists a set $A$ of $\frac{\varphi\left(\frac{p-1}{2}\right)}{\logd_2(1,z)}$ primes, odd and coprime to $z$, such that:

$$\mathcal{G}_S = \left\{ g^{2^x a} m'^{(2^{x+1}-2)a}, (a,x) \in (\{1\} \cup A) \times [\![0, \logd_2(1,z) - 1]\!] \right\}$$

If $p \in \mathcal{P}_{3,4}$, $\mathcal{G}_Z = \left\{ g^{2^x a} m'^{(2^{x+1}-1)a}, (a,x) \in (\{1\} \cup A) \times [\![0, \logd_2(1,z) - 1]\!] \right\}$.

If $p \in \mathcal{P}_{1,4}$, $\mathcal{G}_Z = \left\{ \pm g^{2^x a} m'^{(2^{x+1}-1)a}, (a,x) \in (\{1\} \cup A) \times [\![0, \logd_2(1,z) - 1]\!] \right\}$.

*Proof*: $A$ can be constructed by induction as in remark 2.2.2, by ensuring that $U_0^{(a)}$ is a new element of $\mathcal{G}_S$. Using the identification of $\mathbb{F}_p^*$ to $(\mathbb{Z}/2^n\mathbb{Z}) \times (\mathbb{Z}/z\mathbb{Z})$ and $a, b \in (\{1\} \cup A)$ distinct, suppose $(2ak_0, 2^x al_0) = (2bk_0, 2^y bl_0)$. This means that $2^x al_0 \equiv 2^y bl_0 \ [z]$ thus, multiplying by $2^{\logd_2(1,z)-y}$, $U_0^{(b)}$ is one of the $\left(U_x^{(a)}\right)$, which is impossible.

If $p \in \mathcal{P}_{3,4}$, remember that in this case $\mathcal{G}_S$ and $\mathcal{G}_Z$ have the same size. Thus, once we have $\mathcal{G}_S$, we also get $\mathcal{G}_Z$ by multiplying all elements by the suitable power of $m'$, which is $-1$ anyway.

When $p \in \mathcal{P}_{1,4}$, let us show that elements $\pm g^{2^x a} m'^{(2^{x+1}-1)a}$ are pairwise distinct. Suppose that $m'^a U_x^{(a)} = \pm m'^b U_y^{(b)}$. Then, taking the square yields that $U_{x+1}^{(a)} = U_{y+1}^{(b)}$, which again is impossible.

**Remark 2.2.5**: When $z$ is a prime number, all non-residues $g$ not in $\mathcal{R}^{(n)}$ are generators. If, additionally, $p \in \mathcal{P}_{3,4}$ and $\logd_2(1,z) = z-1 = \varphi(p-1)$ then $A = \emptyset$ for both remark 2.2.2 and proposition 2.2.11.

**Remark 2.2.6**: When $p > 3$ and $z = 1$, any of the $2^{n-1}$ non-residues is a generator, and the sequences $\left(U_x^{(a)}\right)$ are constants. Therefore we need a set $A$ of size $2^{n-2} - 1$ to enumerate all the generators using Proposition 2.2.11. For example, for $p = 257$, we need 63 primes, which means that we must choose primes greater than $p$.

**Remark 2.2.7**: For each generator found, proposition 2.2.9 also gives potential new generators, which can reduce the number of primes needed in remark 2.2.2 or proposition 2.2.11. This approach is left to the reader as it seems less predictable but could prove interesting for values of $p$ greater than those we have tested and such that $\logd_2(1,z)$ is small, since the density of primes goes asymptotically towards 0. Alternatively, we could "sieve" the set $\{k \in [\![1, p-2]\!] \mid \gcd(k, p-1) = 1\}$ by eliminating multiples of prime factors of $p-1$.

### c. Algorithm to obtain the set of primitive roots

The general idea of the algorithm is to search for a semi-primitive root $g \in \mathcal{G}_S$ and the associated $r$ in $\mathcal{R}_p^{(n)} \setminus \mathcal{R}_p^{(n-1)}$ such that $g^z r^{2z} \equiv 1$. Following the proof of proposition 2.2.11, we then construct a set $A$ of primes to obtain the full set $\mathcal{G}_Z$. Let us describe the algorithm a bit more in details below.

**Description of algorithm 1**: We still write $p - 1 = 2^n z$ with $z$ odd. $Q \coloneqq \mathcal{D}_p(\{z\})$, the set of prime divisors of $z$, is also assumed to be known, allowing to calculate $\varphi(z)$.

1. We iterate over the elements $m$ in $[\![2, p - 2]\!]$. For each iteration, two tests are performed, and we go to the next step as soon as one succeeds:
    a. If $m \in \mathcal{G}_{QR}$ and $\forall q \in Q \setminus \{2\}$ $m^{\frac{\varphi(p)}{2q}} \not\equiv 1$ then (Proposition 2.2.6) $m \in \mathcal{G}'_Z$ and we get an element $r \in \mathcal{G}_{QNR}$ (possibly a value $r < m$ on which we have already performed the test). We then set $g_z \equiv m \cdot r^z$ and $g \equiv g_z{}^2$.
    b. Otherwise $m \in \mathcal{G}_{QNR}$ and if, in addition, $\forall q \in Q$ $m^{\frac{\varphi(p)}{q}} \not\equiv 1$ then (Proposition 2.1) $m \in \mathcal{G}_Z$ and we set $g_z = m$ and $g = g_z{}^2$.
2. We iterate over $\mathcal{R}_p^{(n)} \setminus \mathcal{R}_p^{(n-1)} = \{\pm (g_z{}^z)^{2k+1}, 0 \leq k < 2^{n-1}\}$, searching for a value $m'$ of $\mathcal{R}_p^{(n)} \setminus \mathcal{R}_p^{(n-1)}$ such that $g^z m'^{2z} \equiv 1$.
3. We determine $\mathrm{logd}_2(1, z)$ by calculating the first terms of the sequence $(U_x)$ (proposition 2.2.6). The first value of $x \geq 1$ such that $U_x \equiv g$ is indeed $x = \mathrm{logd}_2(1, z)$.
4. We generate odd prime numbers $a_1 \dots a_N \notin \mathcal{D}_p(\{z\})$ such that $U_0^{(a_K)} \notin \{U_x^{(a_k)}, x \in \mathbb{N}, k \in [\![0, K - 1]\!]\}$ (where by convention $U_x^{(a_0)} \coloneqq U_x$) and $N = \left(\varphi\left(\frac{p-1}{2}\right) / \mathrm{logd}_2(1, z)\right) - 1$. For example, these numbers can be obtained from Atkin's algorithm [6] by rejecting those that are not suitable.
5. We have thus recursively constructed $\mathcal{G}_S = \{U_x^{(a_k)}, x \in \mathbb{N}, k \in [\![0, K - 1]\!]\}$ and $\mathcal{G}_Z = \{m'^{a_k} U_x^{(a_k)}\}$ when $p \in \mathcal{P}_{1,4}$, $\mathcal{G}_Z = \{\pm m'^{a_k} U_x^{(a_k)}\}$ otherwise.

## 3. Irreducible quadratic forms

In [2], we studied the density of primes in $E_c = \{X^2 + c, X \in 2\mathbb{N} + r\}$, with $c \in \mathbb{N}^*$ and $r \in \{0,1\}$ such that $r \equiv 1 - c \ [2]$. We empirically corroborated Shanks' conjecture, which gives the asymptotic density of primes in $E_c$:

$$d_{\mathbb{P}|E_c}(x) \sim \frac{h_c}{2\ln(x)}$$

with $h_c = \prod_{p \in \mathbb{P}} \frac{p - t_p(c)}{p - 1} < \infty$ and $t_p(c) = \begin{cases} 0 \text{ if } p \notin \mathcal{D}_p(E_c) \\ 1 \text{ if } p | c \\ 2 \text{ otherwise.} \end{cases}$

Here, we extend the definition of $h_c$ to integer quadratic forms by letting, for $Q = aX^2 + bX + c$:

$$t_p(Q) = |\{x \in \mathbb{F}_p | Q(x) \equiv 0[p]\}|, \quad h_Q = \prod_{p \in \mathbb{P}\setminus\{2\}} \frac{p - t_p(Q)}{p - 1}$$

This definition extends that of $h_c$ i.e. $h_c = h_{(2X+r)^2+c}$.

We continue this study here and present several results about $h_Q$.

## 3.1. Application of Fermat's sum of two squares theorem

**Theorem 2**: Let $c$ be the square of a non-zero integer. We have:

$$\mathcal{D}_p(E_c) = \mathcal{P}_{1,4} \cup \mathcal{D}_p(\{c\}) = \mathcal{P}_{1,4} \sqcup \left(\mathcal{P}_{3,4} \cap \mathcal{D}_p(\{c\})\right)$$

*Proof*: Let us write $c = y^2$. Let $p \in \mathcal{D}_p(E_c)$, and $x \in \mathbb{N}$ such that $p$ divides $x^2 + y^2$. We deduce that $p$ divides $y$ or that $-1$ is a square modulo $p$. But $-1$ is a square modulo $p$ if and only if $p \in \mathcal{P}_{1,4}$, thus we have $p \in \mathcal{P}_{1,4} \cup \mathcal{D}_p(\{y\}) = \mathcal{P}_{1,4} \sqcup \left(\mathcal{P}_{3,4} \cap \mathcal{D}_p(\{y\})\right)$.

Conversely, let $p \in \mathcal{P}_{1,4} \cup \mathcal{D}_p(\{y\})$. If $p \in \mathcal{D}_p(\{y\})$ it is clear that $p$ divides $0^2 + y^2 \in E_c$ so $p \in E_c$. Now assume $p \in \mathcal{P}_{1,4}$. Fermat's sum of two squares theorem [3] ensures that there exist $a, b \in \mathbb{N}^*$ coprime such that $p = a^2 + b^2$. But then we note that $p(c^2 + d^2) = (ac - bd)^2 + (ad + bc)^2$. Since $a$ and $b$ are coprime, there exist $c, d \in \mathbb{Z}$ such that $ad + bc = y$, thus by letting $x = |ac - bd|$ we obtain that $x^2 + y^2$ is a multiple of $p$, hence $p \in E_c$.

**Proposition 3.1a**: If $c$ is a non-zero square, we have:

$$h_c = h_1 \left(\prod_{p \in \mathcal{P}_{3,4} \cap \mathcal{D}_p(\{c\})} \frac{p-1}{p}\right) \left(\prod_{p \in \mathcal{P}_{1,4} \cap \mathcal{D}_p(\{c\})} \frac{p-1}{p-2}\right)$$

*Proof*: We recall that $h_c = \prod_{p \in \mathbb{P}\setminus\{2\}} \frac{p-t_p(c)}{p-1}$. Thus by theorem 2,

$$h_1 = \lim_{n \to \infty} \left(\prod_{\substack{p \in \mathcal{P}_{3,4} \\ p \leq n}} \frac{p}{p-1}\right) \times \left(\prod_{\substack{p \in \mathcal{P}_{1,4} \\ p \leq n}} \frac{p-2}{p-1}\right)$$

and

$$h_c = \lim_{n \to \infty} \left(\prod_{\substack{p \in \mathcal{P}_{3,4} \setminus \mathcal{D}(c) \\ p \leq n}} \frac{p}{p-1}\right) \times \left(\prod_{\substack{p \in \mathcal{P}_{1,4} \setminus \mathcal{D}(c) \\ p \leq n}} \frac{p-2}{p-1}\right).$$

which shows indeed that $\frac{h_c}{h_1} = \left(\prod_{p \in \mathcal{P}_{3,4} \cap \mathcal{D}_p(\{c\})} \frac{p-1}{p}\right) \left(\prod_{p \in \mathcal{P}_{1,4} \cap \mathcal{D}_p(\{c\})} \frac{p-1}{p-2}\right)$.

**Corollary 3.1a**: Let $h \neq h_1 \in \mathbb{R}_+$. If there is a square $c$ such that $h_c = h$ then there exists an infinity of such $c$.

*Proof*: Since the formula in proposition 3.1a depends only on $\mathcal{D}_p(\{c\}) \neq \emptyset$, it suffices to observe that there is an infinity of $c$ which share the same set of prime factors.

**Remark**: If we choose all the prime factors of $c$ in $\mathcal{P}_{1,4}$, then $h_c \geq h_1$.

**Proposition 3.1b**: Let $a, b \in \mathbb{N}$ such that $\mathcal{D}_p(\{a\}) \subset \mathcal{D}_p(\{c\})$ and $b \equiv r\ [2]$. Let $Q = (2aX + b)^2 + c$ and $E_Q = \{Q(x), x \in \mathbb{N}\}$. We assume that $\mathcal{D}_p(\{c\}) \cap \mathcal{D}_p(E_Q) = \emptyset$. Then:

$$h_Q = \lim_{n \to \infty} \left( \prod_{\substack{p \in \mathbb{P} \setminus \mathcal{D}_p(E_Q) \\ p \leq n}} \frac{p}{p-1} \right) \times \left( \prod_{\substack{p \in \mathcal{D}_p(E_Q) \\ p \leq n}} \frac{p-2}{p-1} \right)$$

*Proof*: We must show that $t_p(Q) = 2$ when $p \in \mathcal{D}_p(E_Q)$. As by construction $E_Q \subset E_c$ and $\mathcal{D}_p(\{c\}) \cap \mathcal{D}_p(E_Q) = \emptyset$, $p \in \mathcal{D}_p(E_Q)$ implies $t_p(c) = 2$. Furthermore, $\mathcal{D}_p(\{a\}) \subset \mathcal{D}_p(\{c\})$ implies that $2a$ is invertible modulo $p$, so $t_p(Q) = 2$ too.

**Corollary 3.1b**: Let $c \geq 2$ a square such that $\mathcal{P}_{1,4} \cap \mathcal{D}_p(\{c\}) = \emptyset$. Let $F = \mathcal{D}_p(\{c\})$ and $\alpha \in (\mathbb{N}^*)^F$. We set $p_F^\alpha = \prod_{p \in F} p^{\alpha_p}$, then for any $b \in [\![1, p_F^\alpha - 1]\!]$ coprime to elements of $F$ and of opposite parity to $c$, the quadratic form:

$$Q_{c,\alpha,b} = (2p_F^\alpha X + b)^2 + c$$

verifies $h_{Q_{c,\alpha,b}} = h_1$. Moreover if $b \neq b'$ then the terms of $Q_{c,\alpha,b}$ are disjoint from those of $Q_{c,\alpha,b'}$.

*Proof*: Let $Q = Q_{c,\alpha,b}$. Proposition 3.1a and the assumption $\mathcal{P}_{1,4} \cap \mathcal{D}_p(\{c\}) = \emptyset$ yield $\mathcal{D}_p(E_c) \setminus \mathcal{D}_p(\{c\}) = \mathcal{P}_{1,4}$. However, by construction we also have $\mathcal{D}_p(E_c) \setminus \mathcal{D}_p(\{c\}) = \mathcal{D}_p(E_Q)$, thus proposition 3.1b yields $h_Q = h_1$.

Assume that there exists $b' \in [\![1, p_F^\alpha - 1]\!]$ and $x, y \in \mathbb{Z}$ such that $Q(x) = Q_{c,\alpha,b'}(y)$. We deduce $2p_F^\alpha(x \pm y) = -(b \pm b')$ which implies $b = b'$ and $x = y$. The terms of $Q_{c,\alpha,b}$ and $Q_{c,\alpha,b'}$ are hence disjoint.

**Theorem 3**: There are infinitely many quadratic forms $Q$ with disjoint terms and the same value $h_Q$, which can also be assumed to be arbitrarily close to any value $h \geq h_1$.

Two proofs of this theorem are given.

*First proof*: Proposition 3.1a and its corollary yield that, for any subset $F$ of $\mathcal{P}_{1,4}$, there are infinitely many odd squares such that:

$$h_c = h_1 \left( \prod_{p \in F} \frac{p-1}{p-2} \right)$$

If $c = y^2, c' = y'^2$ are such perfect squares then the equality
$$x^2 + y^2 = x'^2 + y'^2$$
implies $(x - x')(x + x') = y'^2 - y^2$ hence necessary $x$ and $x'$ are smaller than $y'^2 - y^2$. Thus we can construct quadratic forms $Q_k = (2X + m_k)^2 + c_k$ such that $h_{Q_k} = h_{c_k} = h_1 \left( \prod_{p \in F} \frac{p-1}{p-2} \right)$.

Moreover the series $\sum_{p \in \mathcal{P}_{1,4}} \frac{1}{p}$ diverge: indeed Chebotarev's theorem states that $d_{\mathbb{P}|4\mathbb{N}+1}(x) \sim \frac{1}{2\ln(x)}$, so $T_n := \sum_{\substack{p \in \mathcal{P}_{1,4} \\ e^n < p \leq e^{n+1}}} \frac{1}{p} \geq e^{-(n+1)}[\lfloor e^{n+1}\rfloor d_{\mathbb{P}|4\mathbb{N}+1}(e^{n+1}) - \lfloor e^n \rfloor d_{\mathbb{P}|4\mathbb{N}+1}(e^n)] \sim \frac{e-1}{2ne}$ is a divergent series, which implies that $\sum_{p \in \mathcal{P}_{1,4}} \frac{1}{p}$ diverges too. Additionally, $\frac{1}{p} \to_{p \to \infty} 0$. Thus, we can choose $F$ such that $\left(\prod_{p \in F} \frac{p-1}{p-2}\right)$ be arbitrarily close to $h \geq h_1$ fixed.

**Remark**: Similarly, if we choose $F \subset \mathcal{P}_{3,4}$ we can also make $h_Q$ arbitrarily close to any $h \in [0, h_1]$.

*Second proof*: As in corollary 3.1b, for any subset $F$ of $\mathbb{P} \setminus \{2\}$ and for any square $c$ such that $\mathcal{D}_p(\{c\}) \subset F$, we can get $a, b$ such that:

$$h_{(2aX+b)^2+c} = h_1 \left(\prod_{p \in F \cap \mathcal{P}_{1,4}} \frac{p-1}{p-2}\right).$$

We conclude as in the first proof: we can either take an infinity of values for $c$ as in the first proof, or build a sequence $(a_k, b_k)$ such that $a_{k+1}$ is a multiple of $a_k$ and for any $l \leq k$, $b_l \not\equiv b_{k+1} [a_l]$.

**Remark**: If $p \in F \cap \mathcal{P}_{1,4} \setminus \mathcal{D}_p(\{c\})$, to "eliminate" it, we need to determine the solutions of $X^2 \equiv -c [p]$. This equation is discussed in the next section.

By analogy with Shanks' conjecture, $h_Q$ can be conjectured to be linked with the density of prime numbers of the form $Q(x)$ as in the following generalisation:

**Conjecture 3.1**: Let $Q$ be an irreducible quadratic form and $E_Q = \{Q(x), x \in \mathbb{N}\}$. Then:

$$d_{\mathbb{P}|E_Q}(x) = \frac{h_Q}{\ln(x)} + o\left(\frac{1}{\ln(x)}\right)$$

where $d_{\mathbb{P}|E_Q}$ is the density of primes in $E_Q$ (see [2] for a formal definition).

### 3.2. Tonelli-Shanks algorithm

There are several algorithms to calculate the square root of an integer modulo $p$, e.g. Tonelli-Shanks [4][5], Cipolla [7] and Daniel Bernstein [8] algorithms. In this section, we rewrite the Tonelli-Shanks algorithm [4][5] to solve equation $(3.2)$ $x^2 + c \equiv 0 [p]$ simultaneously for several values of $c$.

**Proposition 3.2**: As in section 2, let us write $p - 1 = 2^n z$ with $z$ odd. If $m$ is a quadratic residue, there exists $r \in \mathcal{R}^{(n)}$ such that:

$$\left(m^{\frac{p-z}{2}} r\right)^2 \equiv m$$

The square roots of $m$ are thus $\pm m^{\frac{p-z}{2}} r$ and $m^{\frac{p-z}{2}} r$ is a quadratic residue if and only if $r$ also is.

*Proof*: We still identify $\mathbb{F}_p^*$ to $(\mathbb{Z}/2^n\mathbb{Z}) \times (\mathbb{Z}/z\mathbb{Z})$. $m$ is a quadratic residue and therefore corresponds to an element in form $(2k, l)$. We know that $p - z = 1 + (2^n - 1)z$ is even and $p - z \equiv 1\ [z]$ hence $\frac{p-z}{2} \times 2 \equiv 1\ [z]$. On the other hand, modulo $2^n$, $\frac{p-z}{2}.2k \equiv (p-z)k \equiv (1-z)k$. Let $r \in \mathbb{F}_p^*$ which corresponds to $(zk, 0)$. Then $m^{\frac{p-z}{2}}r$ corresponds to $\left(k, \frac{p-z}{2}l\right)$, so its square is so equal to $m$. Moreover, $m^{\frac{p-z}{2}}r$ is a quadratic residue if and only if $k$ is even, which is equivalent to $zk$ being even, i.e. to $r$ being quadratic residue.

**Remark 3.2**: If $p \in \mathcal{P}_{1,4}$, we can also write the square root of $m$ as $m^{\frac{p+1-2z}{4}}r$. Indeed, on one hand $\frac{p+1}{2}$ and $z$ are both odd with $z \leq \frac{p-1}{4} < \frac{p+1}{2}$ so $\frac{p+1-2z}{4} \in \mathbb{N}^*$, and on the other hand:

$$\left(m^{\frac{p+1-2z}{4}}r\right)^2 \equiv m^{\frac{p-1}{2}}m^{p-z}r^2 \equiv 1 \times \left(m^{\frac{p-z}{2}}r\right)^2 \equiv m\ [p]$$

**Corollary 3.2**: If $p \in \mathcal{P}_{3,4}$ and $m \in \mathcal{G}_{QR}$, then its square roots are $\pm m^{\frac{p+1}{4}}$.

*Proof*: This follows directly from proposition 3.2, since then $r = \pm 1$ and $p - z = \frac{p+1}{2}$. We can also verify directly $m^{\frac{p+1}{2}} = m^{\frac{p-1}{2}+1} = m$ by Euler's criterion, since $m \in \mathcal{G}_{QR}$.

Searching for the roots of $X^2 + c$ modulo $p$ allows us to obtain, if it exists, the smallest multiple of $p$ in $E_c$. We describe such an algorithm below.

**Description of algorithm 2**: Search for $(x_1 \ldots x_n)$ such that $x_k$ is the smallest non-negative integer verifying $x_k^2 + c_k \equiv 0\ [p]$ with $x_k$ and $c$ of opposite parity (i.e. $x_k^2 + c_k$ is the smallest multiple of $p$ in $E_{c_k}$), for a $n$-uplet $(c_1 \ldots c_n) \in (\mathbb{N}^*)^n$ and $p \in \mathbb{P} \setminus \{2\}$ fixed. We write $p - 1 = 2^n z$ with $z$ odd. By convention, we let $x_k = \infty$ if $X^2 + c_k$ has no root modulo $p$.

1. If $p \in \mathcal{P}_{3,4}$, for each value of $c_k$, we compute $c_k^{\frac{p-1}{2}}$ modulo $p$. If $c_k^{\frac{p-1}{2}} \equiv 1\ [p]$, we set $x_k = \infty$, otherwise we set $x_k$ the smallest positive integer with parity opposite to that of $c_k$ and such that $x_k \equiv \pm c_k^{\frac{p+1}{4}}\ [p]$.

2. If $p \in \mathcal{P}_{1,4}$, we set $g = 2$ and increase it until $g \in \mathcal{G}_{QNR}$ (i.e. by Euler criterion $g^{\frac{p-1}{2}} \equiv -1$) then:
   a) We enumerate $\mathcal{R}_p^{(n-1)} = \{\pm(g^z)^{2k}[p], 1 \leq k \leq 2^{n-2}\}$.
   b) For each value $c_k$:
      a. If $c_k^{\frac{p-1}{2}} \equiv -1\ [p]$, we set $x_k = \infty$.
      b. Otherwise for each element $r$ of $\mathcal{R}_p^{(n-1)}$, we test if $(-c)^{\frac{p+1-2z}{2}}r^2 \equiv -c\ [p]$ modulo $p$. If true, we set $x_k$ the smallest positive integer with parity opposite to that of $c_k$ and such that $x_k \equiv \pm(-c_k)^{\frac{p+1-2z}{4}}r$.

## 4. CONCLUSION
Based on the study of primitive roots ($\mathcal{G}_Z$) and semi-primitive roots ($\mathcal{G}_S$) modulo $p \in \mathbb{P} \setminus \{2\}$, and in particular thanks to recursive sequences of elements of $\mathcal{G}_Z$ and $\mathcal{G}_S$, we proposed an original algorithm returning these two sets in full, which is based on a prime number generator algorithm rather than on $gcd$ calculations.

We then came back to Shanks' conjecture, which we generalised to irreducible quadratic forms, and showed that the conjectured asymptotic density constant could be well controlled. To this end, we rewrote the Tonelli-Shanks algorithm for solving the equation $X^2 + c \equiv 0 \ [p]$.

Acknowledgments: We would like to thank François-Xavier VILLEMIN for his attentive comments and suggestions.